\documentclass[letterpaper, 10 pt, conference]{ieeeconf}  

\IEEEoverridecommandlockouts                              
\usepackage{longtable}
\usepackage{comment}
\usepackage{amsmath}
\usepackage{graphicx}
\usepackage{color}
\usepackage[dvips]{epsfig}

\usepackage{graphicx, amssymb}
\usepackage{amsmath}

\newcommand{\beq}{\begin{equation}}
\newcommand{\eeq}{\end{equation}}
\newcommand{\carre} {\hfill $\blacksquare$}

\newtheorem{theo}{Theorem}

\newtheorem{lem}{Lemma}
\newtheorem{corr}{Corollary}

\newtheorem{pro}{Proposition}

\title{\LARGE \bf
Controllability Analysis of Threshold Graphs and Cographs}

\author{Shima Sadat Mousavi$^{\dagger}$, Mohammad Haeri$^{\dagger}$, and Mehran Mesbahi$^{\ddagger}$
\thanks{$^{\dagger}$The authors are with the Department of Electrical Engineering, Sharif University of Technology, Tehran, Iran. Emails: shimasadat$\_$mousavi@ee.sharif.edu, haeri@sharif.ir.}%
\thanks{$^{\ddagger}$The author is with the Department of Aeronautics and Astronautics, University of Washington, WA 98195. Email: mesbahi@aa.washington.edu.}%
}
\begin{document}

\maketitle
\thispagestyle{empty}
\pagestyle{empty}

\begin{abstract}
In this paper, we investigate the controllability of a linear time-invariant network following a Laplacian dynamics defined on a threshold graph. In this direction, an algorithm for deriving the modal matrix  associated with the Laplacian matrix for this class of graphs is presented. Then,  based on the Popov-Belevitch-Hautus criteria, a procedure for the selection of control nodes  is proposed. The procedure involves partitioning the nodes of the graph into cells with the same degree; one node from each cell is then selected.  We show that the remaining nodes can be chosen as the control nodes rendering the network controllable. Finally, we consider a wider class of graphs, namely cographs, and examine their controllability properties.  
\end{abstract}


\section{Introduction}

Networks are the backbone of modern society.
Social networks, the internet, and energy networks, are examples of some of the critical networks that we rely on their operation in our daily lives. 
As such, the control, security, and management of these and other types of networks are of paramount importance, providing a rich class of system theoretic questions for the control community~\cite{mesbahi2010graph}.
One foundational class of questions on networked systems
pertain to their controllability~\cite{tanner2004controllability, liu2011controllability,mousavi2017structural}. 
Controllability analysis on networks can also provide a framework for designing network topologies with favorable controllability properties.
However, some of the basic controllability questions on networks-even for the linear case-are nontrivial.
For example, finding a minimum cardinality set of control nodes that ensures the controllability of a large-scale network  through the classical rank conditions is NP-hard. 
Accordingly, an alternative means of examining network controllability is via its topological properties. In this direction,  controllability analysis of networks with the so-called Laplacian dynamics has received a lot of attention, primary due to their role in consensus-type collective behaviors such as synchronization \cite{rahmani2009controllability, zhang2011controllability, zhang2014upper, yaziciouglu2016graph}.  

The results in the literature on the controllability analysis of networks with Laplacian dynamics can be classified into two categories. In the first category,  a general topology has been  considered  for the network, and   necessary or sufficient conditions for its controllability from a graph-theoretic point of view have been  presented. 
These conditions have been stated in terms of notions such as graph symmetry \cite{rahmani2009controllability,chapman2015state}, equitable partitions \cite{rahmani2009controllability, egerstedt2012interacting,martini2010controllability,zhang2011controllability, zhang2014upper, cao2013class,aguilar2017almost}, distance partitions \cite{zhang2011controllability, zhang2014upper}, and pseudo monotonically increasing sequences \cite{yazicioglu2012tight, yaziciouglu2016graph}. However, some of the approaches have a few important limitations. For example, none of these conditions are necessary \emph{and} sufficient for network controllability; rather, they are used in deriving lower or/and upper bounds on the dimension of the controllable subspace. More importantly, these results cannot be utilized for efficient selection of control nodes rendering a network controllable. For example, it is known that the existence of a symmetry in the structure of a network with respect to its control nodes is destructive to its controllability \cite{rahmani2009controllability}, but finding a minimum cardinality set of nodes breaking all symmetries for general networks is NP-hard \cite{chapman2015state}.  

The second category of existing works includes those that consider controllability of special classes of networks~\cite{aguilar2015graph}. For example, controllability of networks with path graphs \cite{parlangeli2012reachability, mousavi2016controllability}, cycle graphs \cite{parlangeli2012reachability}, complete graphs \cite{zhang2011controllability}, circulant graphs \cite{nabi2013controllability}, multi-chain graphs \cite{hsu2017necessary}, grid graphs \cite{notarstefano2013controllability}, and tree graphs \cite{ji2012leaders} have already been explored. 
In these cases, stronger conditions for network controllability can be derived. 
In particular, for some of these graphs, the minimum number of control nodes from which the associated network is controllable has been determined. 
Note that the stronger controllability conditions derived for these special classes of graphs are resulted from a better characterization of the eigenvectors associated with their Laplacian. In fact, based on the Popov-Belevitch-Hautus (PBH) test, the controllability of a system solely depends on its associated eigenvectors and how they relate to the input structure. Subsequently, by identifying the eigenspace of the network (i.e., the space of eigenvectors associated with each eigenvalue of the Laplacian matrix), the controllability problem can be addressed.

Adopting a similar approach, in this paper, we consider the controllability problem for the Laplacian networks defined on cographs. 
Cographs have been independently introduced by different research works, and as such, admit a few equivalent definitions. 
For example, there is no subgraph isomorphic to a path of size four in cographs. 
Moreover, cographs can be generated by successively operating \emph{joins} and \emph{unions} among isolated nodes \cite{biyikoglu2007laplacian}.  
Cographs have many applications in diverse areas of computer science and mathematics \cite{corneil1984cographs}. 
Moreover, they include other known classes of graphs with special structures.
 For example, threshold graphs with applications in areas such like modeling social and psychological networks, synchronizing parallel processes, and cyclic scheduling problems, are cographs~\cite{mahadev1995threshold}. 
 There are different representations for threshold graphs as well; for instance, threshold graphs can be uniquely determined by a binary construction sequence~\cite{hagberg2006designing}. 
 In \cite{aguilar2015laplacian}, the controllability of a threshold graph from only a single control node has been explored. In particular, in this work it has been proved that a threshold graph is controllable from a single controller only if it is an antiregular graph with $n-1$ different degrees.
 Subsequently, the work~\cite{hsu2016controllability} extended the result of \cite{aguilar2015laplacian} by considering threshold graphs with only one repeated degree. 

The main contributions of the present paper are as follows: First,  we consider a very general threshold graph and allow it to have any number of repeated degrees. In this regard, we assume that  a threshold graph is  described by its construction sequence and derive a modal matrix associated with its Laplacian. Then, we explore the controllability of a network defined on this graph.
By adopting an approach different from the one  used in \cite{hsu2016controllability}, we show that for any repeated degree, by independently controlling any node of that degree except one, we can ensure controllability of the network. In particular, we prove that the minimum number of control nodes to fully control a threshold graph is the difference between the size of the network and the number of its distinct degrees. Moreover, we present a systematic method to choose the control nodes.  Next, we provide a controllability analysis of a cograph via its eigenspace. In this direction, we provide a method for deriving an  input matrix with the minimum rank that renders the network controllable. 

The organization of the paper is as follows. First, the notation and preliminaries are provided. In \S III, the eigenvectors associated with a threshold graph are derived, and necessary and sufficient conditions for the controllability of this class of graphs is established. %
\S IV is dedicated to the controllability analysis of networks on cographs. Finally, \S V concludes the paper.

\section{Notation and Preliminaries}
In this section, the notation and preliminaries for our 
subsequent discussion is presented.\\%
\newline
\emph{Notation:} The set of real numbers and integer numbers are respectively, denoted by $\mathbb{R}$ and $\mathbb{Z}$. For a matrix $M\in\mathbb{R}^{p\times q}$, $M_{ij}$ is the entry of $M$ in its $i$th row and $j$th column. Furthermore, $M_{i,:}$ and $M_{:,j}$ represent the $i$th row and $j$th column of $M$.
The $n\times n$ identity matrix is denoted by $I_n$, and $e_j$ represents its $j$th column. The vectors of all 1's  and all 0's with size $n$ are respectively, denoted by  $\textbf{1}_n$ and $\textbf{0}_n$. Also, an $n\times m$ matrix of all 1's (resp., 0's) is given by $\textbf{1}_{n\times m}$ (resp., $\textbf{0}_{n\times m}$).
For a set $\mathcal{S}$, we denote its cardinality by $|\mathcal{S}|$. 


\textit{Graph:} A graph\footnote{All graphs in this paper are assumed to be undirected, unweighted, and loop-free.} $G$ of size $n$ is represented by $G=(V,E)$, where 
$V=\{1,\ldots,n\}$ is its node set, and $E$ denotes its edge set. 
The node $j$ is called a neighbor of the node $i$ if $\{i,j\}\in E$. We denote by $N(i)$ the set of neighbors of  $i\in V$. The degree of the node $i$ is defined as $d(i)=|N(i)|$. The degree matrix of the graph $G$ is defined as $\Delta(G)=\mathrm{diag}(d(1),\ldots,d(n))$. Then, the Laplacian matrix $L(G)$ is given by $L(G)=\Delta(G)-A(G)$, where $A(G)$ is the (0,1)-adjacency matrix associated with  the graph $G$. 
The degree sequence $D(G)$  is a nondecreasing sequence of node degrees of $G$, which is  defined as $D(G)=(d_1,\ldots,d_n)$, where $d_1\leq \ldots \leq d_n$. 
Let $s$ be the number of distinct degrees in $D(G)$. Then, one can write $D(G)=(\tilde{d}_1^{p_1},\ldots\tilde{d}_s^{p_s})$, where $\tilde{d}_1\leq \ldots\leq \tilde{d}_s$ are the $s$ distinct degrees of the nodes, and $p_i$ is the multiplicity of the degree $\tilde{d}_i$, $1\leq i\leq s$, among the nodes of $G$.
 
 \emph{Eigenpairs:} With a slight abuse of notation, by eigenvalues and eigenvectors of a graph $G$, we mean the eigenvalues and eigenvectors of its Laplacian matrix $L(G)$. Since $L(G)\geq0$, all of its eigenvalues are real and nonnegative. Let $\Lambda(G)=(\lambda_1,\ldots,\lambda_n)$ be the spectrum of the graph $G$, where $\lambda_1\leq \lambda_2\leq\ldots\lambda_n$. Then, $\lambda_1=0$, and if $G$ is connected, we have $\lambda_2\neq 0$. If $\tilde{\lambda}_1\leq\ldots\leq \tilde{\lambda}_r$ are the $r$ distinct nonzero eigenvalues of $G$, then for a connected $G$, we can write  $\Lambda(G)=(0,\tilde{\lambda}_1^{q_1},\ldots,\tilde{\lambda}_r^{q_r})$, where $q_i$, $1\leq i\leq r$, is the algebraic multiplicity of the eigenvalue $\tilde{\lambda}_i$. Then, $\mathcal{M}=\mbox{max}\{q_1,\ldots, q_r\}$ is the maximum multiplicity of eigenvalues of $G$. We denote an eigenpair of the graph $G$ by the pair $(\lambda_i,\nu_i)$, where $L(G)\nu_i=\lambda_i\nu_i$, $1\leq i\leq n$.  The vector $\nu_i\in \mathbb{R}^n$ is an eigenvector of $G$  associated with the eigenvalue $\lambda_i$. One can see that every graph has $(0,\textbf{1}_n)$ as one of its eigenpairs. Now, assume that $\tilde{\lambda}_i$ is an eigenvalue of $G$ with multiplicity $q_i$. Then, there are $q_i$ independent eigenvectors $\nu^{(i)}_1,\ldots,\nu^{(i)}_{q_i}$ associated with $\tilde{\lambda}_i$. Let $V^{(i)}=[\nu^{(i)}_1,\ldots,\nu^{(i)}_{q_i}]\in\mathbb{R}^{n\times q_i}$. Then every eigenvector $\nu$ associated with $\tilde{\lambda}_i$ can be written as $\nu=V^{(i)}C$, for some $C\in\mathbb{R}^{q_i}$. Let the nonsingular matrix $V(G)=[\nu_1,\ldots,\nu_n]\in\mathbb{R}^{n\times n}$ be a \emph{modal} matrix associated with $\Lambda(G)$, where $L(G)V(G)=V(G)\mbox{diag}(\Lambda(G))$. 
We can also consider an unordered sequence of eigenvalues of $L(G)$. Let $\bar{\Lambda}(G)=(\bar{\lambda}_1,\ldots,\bar{\lambda}_n)$ be a sequence of eigenvalues of $L(G)$, not necessarily ordered in a nonincreasing or nondecreasing way. Moreover, we define $\bar{V}(G)=[\bar{\nu}_1,\ldots,\bar{\nu}_n]$ as a modal matrix associated with $\bar{\Lambda}(G)$, where $(\bar{\lambda}_i,\bar{\nu}_i)$, $1\leq i\leq n$, is an eigenpair of $G$. 

\subsection{Cographs and Threshold Graphs}

We now introduce the notion of cographs and threshold graphs; we also provide theorems about their corresponding spectrum.

Let $G_1=(V_1 ,E_1)$ and $G_2=(V_2,E_2)$ be two disjoint graphs of respectively, sizes $n_1$ and $n_2$. The \emph{union} of the two graphs is a graph of size $n=n_1+n_2$, which is defined as $G_1+G_2=(V_1\cup V_2, E_1\cup E_2)$. Moreover, the \emph{join} of the two graphs represented by $G_1 * G_2$ is obtained from $G_1+G_2$ by adding new edges from each node of $G_1$ to any node of $G_2$. A graph is called a \emph{cograph} (or a \emph{decomposable graph}) if it can be constructed from isolated nodes by successively performing the join and union operations. 

Now, let $\bar{\Lambda}(G_{1})=(0, \bar{\psi}_2,\ldots,\bar{\psi}_{n_1})=(0,\bar{\Psi})$ and $\bar{\Lambda}(G_{2})=(0, \bar{\pi}_2,\ldots,\bar{\pi}_{n_2})=(0,\bar{\Pi})$, with $\bar{\Psi}=( \bar{\psi}_2,\ldots,\bar{\psi}_{n_1})$ and $\bar{\Pi}=(\bar{\pi}_2,\ldots,\bar{\pi}_{n_2})$, be unordered sequences of eigenvalues of $G_1$ and $G_2$.  Moreover, let  $\bar{V}(G_{1})=(\textbf{1}_{n_1}, \bar{u}_2,\ldots,\bar{u}_n)=(\textbf{1}_{n_1},\bar{U})$ and $\bar{V}(G_{2})=(\textbf{1}_{n_2}, \bar{w}_2,\ldots,\bar{w}_n)=(\textbf{1}_{n_2},\bar{W})$ be respectively the modal matrices of $G_1$ and $G_2$ associated with $\bar{\Lambda}(G_1)$ and $\bar{\Lambda}(G_2)$. Then by the next result, one can establish the eigenvalues and eigenvectors of the join and the union of $G_1$ and $G_2$. Consider a vector $v\in \mathbb{R}^n$ and a scalar $m\in \mathbb{R}$. With a slight abuse of notation in this theorem,  we let $v+m=v+m\textbf{1}_{n}$.

\begin{theo}[\cite{merris1998laplacian}]
For two graphs $G_1$ and $G_2$ of respectively, sizes $n_1$ and $n_2$, we have:
\begin{equation*}L(G_1+G_2)=\begin{bmatrix}L(G_1) & \textbf{0}_{n_1\times n_2}\\ \textbf{0}_{n_2\times n_1} & L(G_2)\end{bmatrix},\end{equation*} 
 \begin{equation*}L(G_1*G_2)=\begin{bmatrix}L(G_1)-n_2I_{n_1} & \textbf{0}_{n_1\times n_2}\\ \textbf{0}_{n_2\times n_1} & L(G_2)-n_1I_{n_2}\end{bmatrix},\end{equation*}
\begin{equation*}\bar{\Lambda}(G_1+G_2)=(0, \bar{\Psi}, \bar{\Pi}, 0),\end{equation*}
\begin{equation*}\bar{\Lambda}(G_1*G_2)=(0, \bar{\Psi}+n_2, \bar{\Pi}+n_1, n+m),\end{equation*}
\begin{equation*}\bar{V}(G_1+G_2)=\begin{bmatrix}\textbf{1}_{n_1} & \bar{U} & \textbf{0}_{n_1} & n_2\textbf{1}_{n_1}\\\textbf{1}_{n_2} & \textbf{0}_{n_2} & \bar{W}& -n_1\textbf{1}_{n_2}\end{bmatrix},\end{equation*}
\begin{equation*}\bar{V}(G_1*G_2)=\bar{V}(G_1+G_2).\end{equation*}
\label{th3}
\end{theo}

 Now, let us start with one isolated node as the initial graph, and in each step, connect an isolated node to the former graph through the join or union operation. The obtained graph is referred to as a \emph{threshold graph}, which is a spacial type of a cograph. One can associate a binary \emph{construction sequence} $T^G\in\{0,1\}^n$ to a threshold graph $G$ of size $n$, where $T^G(1)=0$, and for $1<i\leq n$, $T^G(i)=0$ (resp., $T^G(i)=1$) if the node $i$ is added to the former graph by the union (resp. join) operation \cite{hagberg2006designing}. In fact, any threshold graph $G$ can be uniquely determined by its construction sequence $T^G$. In this paper, we assume that all threshold graphs are described and given by their associated construction sequences. Adopting this notation, the next theorem connects the spectrum of a threshold graph to its degree sequence. 
 
\begin{theo}[\cite{hammer1996laplacian}]
In a threshold graph with the the degree sequence $D(G)$, $\lambda_{n-i+1}=|\{j:d(j)\geq i\}|$. Now, let 
$D(G)=(\tilde{d}_1^{p_1},\ldots\tilde{d}_s^{p_s})$ and $\Lambda(G)=(0,\tilde{\lambda}_1^{q_1},\ldots,\tilde{\lambda}_s^{q_s})$. Then if $s=2l$ (resp., $s=2l+1$), for some $l\in \mathbb{Z}$, we have: \\
\begin{align*} \tilde{\lambda}_i=\begin{cases} \tilde{d}_i, \quad 1 \leq i\leq l,\\ \tilde{d}_i+1, \quad l+1\leq i\leq s \end{cases},\end{align*}
\begin{align*}  q_i=\begin{cases} p_i-1, \quad i=l,\ (\mbox{resp.,} \ i=l+1)\\p_i \quad \mbox{otherwise}\end{cases}.\end{align*}
\label{th2}
\end{theo}
  
\subsection{Problem Formulation}

In this paper, we consider a linear time-invariant (LTI) network with the graph structure $G$ and the so-called Laplacian dynamics described as:

\begin{equation}
\dot{x}=Ax+Bu,
\label{e1}
\end{equation} 
where $A=-L(G)$, and $L(G)\in \mathbb{R}^{n\times n}$ is the Laplacian matrix associated with the graph $G$. Moreover, $x=[x_1,\ldots,x_n ]^T$ is the vector of states of the nodes, and $u=[u_1,\ldots,u_m ]^T$ is the vector of input signals. Also, $B\in \mathbb{R}^{n\times m}$ is the input matrix whose nonzero entries determine the nodes where the input signals are directly injected. In this paper, $G$
 is assumed to be  a threshold graph or a general cograph, and the controllability of the network is investigated. 
Specifically, we provide conditions ensuring the controllability of the network  and find the minimum number of independent input signals (or controllers) that render the network controllable. 
 
In the first step, we assume that any input signal can be injected into only one node, referred to as the control node. Thus, the input matrix $B$ can be defined as 

\begin{equation}
B=[e_{j_1},\ldots,e_{j_m}],
\label{B}
\end{equation}
 where $j_i\in\{1,\ldots,n\}$, for $1\leq i\leq m$, and $V_C=\{j_1,\ldots,j_m\}$ is the set of control nodes in the network. In the next step, we consider a general matrix $B$ whose entries can be any real numbers. Then, we find an input matrix $B$ with the minimum number of columns (or equivalently, the minimum number of independent inputs) that renders the network controllable.

In order to investigate the controllability of networks, we use the PBH controllability test as follows.
 
 \begin{pro}[\cite{sontag2013mathematical}]
 A system with dynamics (\ref{e1}) (or the pair $(A,B)$)  is controllable if and only if for any nonzero (left) eigenvector $\nu$ of  $A$, we have $\nu^TB\neq 0$. 
 \end{pro}
 
 The PBH test can be stated in another equivalent way: a system with dynamics (\ref{e1}) is controllable if and only if for every eigenvalue $\lambda$ of $A$, the matrix $[\lambda I_n-A, B]\in \mathbb{R}^{n\times (n+m)}$ is full rank.
 
Note that if we want to select the set of control nodes for a network of size $n$ by relying on the PBH test, we are required to adopt a brute-force 
exponential time algorithm, that is
computationally infeasible for large-scale networks. 
In this paper, we provide controllability conditions for special classes of graphs that can be efficiently inferred from the corresponding network topology.


\section{Controllability of Networks with Threshold Graphs}

In this section, we investigate the controllability of networks with dynamics (\ref{e1}) whose structures are described by threshold graphs. To this aim, the eigenspace of a threshold graph are first examined.
  
\subsection{Eigenspace of a Threshold Graph}

We consider a construction sequence $T^G$ associated with a threshold graph $G$ and  proceed to characterize the eigenvalues and eigenvectors of its Laplacian matrix.

As mentioned previously, considering the sequence $T^G$, a threshold graph $G$ of size $n$ can be constructed in $n$ steps, where in each step, an isolated node is added to the graph through the join or union operation. Consider the node added to $G$ in the $i$th step, and let it be indexed as $i$, $1\leq i\leq n$. 
By the next result, given  $T^G$, one can provide the node degrees of $G$. 

\begin{pro}
Consider the construction sequence $T^G$ associated with a threshold graph $G$. Then, for every $1\leq i\leq n$, $d(i)=T^G(i)\times(i-1)+|\{i<j\leq n: T^G(j)=1\}|.$
\label{pr2}
\end{pro}

\emph{Proof:} For some $1\leq i\leq n$, first let $T^G(i)=1$. Then node $i$ is added to the set of nodes $\{1,\ldots,i-1\}$ through the join operation. In other words, it is connected to all nodes $j$ where $1\leq j<i$. Moreover, for a node $k$ where $k>j$, if  $T^G(k)=1$, $\{i,k\}\in E$, and if  $T^G(k)=0$, then $\{i,k\}\notin E(G)$. Thus, $d(i)=(i-1)+|\{i<j\leq n: T^G(j)=1\}|$. On the other hand, if $T^G(i)=0$, the node $i$ is connected only to the nodes which are added to the graph through a join operation in some step $j$, where $j>i$. In other words, $\{i,k\}\in E$ if $k>i$ and $T^G(k)=1$, which completes the proof. 
\carre

From Proposition \ref{pr2}, with a construction sequence, we can find the degree sequence of the associated threshold graph. Moreover, the next result provides conditions on elements of $T^G$ under which two nodes $i$ and $j$  have the same degree. 

\begin{lem}
Consider the construction sequence $T^G$ associated with a threshold graph $G$. For some nodes $i,j\in V$, where $1\leq i<j\leq n$, we have $d(i)=d(j)$ if and only if one of the following three conditions holds:
\begin{enumerate}
\item $(T^G(i), T^G(i+1),\ldots,T^G(j))=(0,0,\ldots,0)$,\\
\item $(T^G(i), T^G(i+1),\ldots,T^G(j))=(1,1,\ldots,1)$,\\
\item $i=1$, and $(T^G(1), T^G(2),\ldots,T^G(j))=\\(0,1,\ldots,1)$.
\end{enumerate}
\label{lem1}
\end{lem}

\emph{Proof:} For the sufficiency part, using Proposition \ref{pr2}, one can verify that if any of the three conditions hold, $d(i)=d(j)$. Now, let us prove the necessity by contradiction. Let $d(i)=d(j)$. First, assume that $T^G(i)=T^G(j)$, but for some $i<k<j$, $T^G(i)\neq T^G(k)$. Then, if $T^G(i)=0$, by Proposition \ref{pr2}, $d(i)>d(j)$, and if $T^G(i)=1$, $d(j)>d(i)$, which contradicts the assumption. Now, assume that $T^G(i)=1$ and $T^G(j)=0$. Then, $d(i)\geq (i-1)+d(j)$. Moreover, since $T^G(1)=0$, we have $i>1$. Then, $d(i)>d(j)$. On the other hand, let $i>1$ and $T^G(i)=0$, while $T^G(j)=1$. Define $k_1=|\{i<k<j:k=1\}|$ and $k_2=|\{j<k\leq n:k=1\}|$. Then, $d(i)=k_1+k_2+1$ and $d(j)=j-1+k_2\geq k_1+k_2+i$. Then, since $i>1$, $d(j)>d(i)$, which is a contradiction. 
\carre

By Lemma \ref{lem1}, we can also conclude that the degrees of the nodes 1 and 2 in a threshold graph are the same.
\begin{corr}
In a threshold graph $G$ whose nodes are ordered and indexed based on the construction sequence $T^G$, we have $d(1)=d(2)$.
\label{corr1}
\end{corr}

\emph{Proof:} Since $T^G(1)=0$,  then $(T^G(1), T^G(2))=(0,1)$ or $(0,0)$. Accordingly, from Lemma \ref{lem1}, $d(1)=d(2)$.
\carre

By applying Proposition \ref{pr2}, for a threshold graph $G$ with a construction sequence $T^G$, one can obtain the degree sequence $D(G)$. Then, based on Theorem \ref{th2}, the ordered nondecreasing sequence of eigenvalues of $G$ denoted by $\Lambda(G)$ is provided. In Fig. \ref{f1}, an algorithm that generates the modal matrix $V(G)$ associated with $\Lambda(G)$ is presented. Let $\mathcal{C}:\{1,\ldots,n\}\rightarrow \mathbb{R}^{n}$ such that $\mathcal{C}(i)=\begin{bmatrix}\textbf{1}_i^T& -i & \textbf{0}^T_{n-i-1} \end{bmatrix}^T$.

 \begin{figure}[htb]
\mbox{}\hrulefill
\vspace{-.4em}
\\
\textbf{Algorithm 1:}

\vspace{-.7em}
\mbox{}\hrulefill
\vspace{-.01em}
\\
\small \textbf{Input:} The construction sequence $T^G$\\
\textbf{Output:} The modal matrix $V(G)$\\
$V_{:,1}(G)=\textbf{1}_n$\\
$k=1$\\
$r=0$\\
\begin{tabular}{l}
\textbf{for} $i=n:-1:2$\\
\begin{tabular}{|l}
�\textbf{if} $T^G(i)=0$\\
\begin{tabular}{|l}
$k=k+1$\\
$V_{:,k}(G)=\mathcal{C}(i-1)$\\
\end{tabular}
\\
\textbf{else}\\
\begin{tabular}{|l}
$V_{:,n-r}=\mathcal{C}(i-1)$\\
$r=r+1$\\
\end{tabular}
\\
\textbf{end if}\\
\end{tabular}
\\
\textbf{end for}\\
\end{tabular}\\
\textbf{return} $V(G)$\\
\vspace{-.5em}
\mbox{}\hrulefill
\caption{An algorithm that generates the modal matrix $V(G)$ associated with the construction sequence $T^G$ of a threshold graph $G$.}
\vspace{-0.08in}
\label{f1}
\end{figure}

\begin{theo}
For a threshold graph $G$ with a given construction sequence $T^G$, the matrix $V(G)$ obtained by Algorithm 1 is the modal matrix of $G$ associated with the spectrum $\Lambda(G)=(0,\lambda_2,\ldots, \lambda_n)$ (s.t., $0\leq \lambda_2\leq \ldots \leq \lambda_n$).
\label{th4}
\end{theo}

\begin{figure}[hbt]
\includegraphics[width=.23\textwidth]{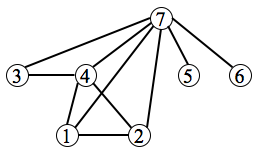}
\centering
\caption{A threshold graph $G$ associated with $T^G=(0,1,0,1,0, 0,1)$.}
\label{ex1}
\end{figure}

Before presenting the  proof of  Theorem \ref{th4}, let us consider a sample run of Algorithm 1 for a construction sequence $T^G=(0,1,0,1,0, 0,1)$ associated with a threshold graph $G$. The graph $G$ is shown in Fig. \ref{ex1}. The nodes are indexed according to the number of the step in which they are added to the graph through the join or union operation. Using Proposition \ref{pr2}, one can find that $d(1)=3$, $d(2)=3$, $d(3)=2$, $d(4)=4$, $d(5)=1$, $d(6)=1$, and $d(7)=6$.  Moreover, from Theorem \ref{th2}, the spectrum of the graph is obtained as $\Lambda(G)=(0,1,1,2,4,5,7)$. 
Now, let us run Algorithm 1 to generate $V(G)$. The first column of $V(G)$ is $\textbf{1}_7$. One can see that for $i=6,5,3$, $T^G(i)=0$. Then, the 2nd, the 3rd, and the 4th columns of $V(G)$ are respectively, equal to $\mathcal{C}(5)=(1,1,1,1,1,-5,0)^T$, $\mathcal{C}(4)=(1,1,1,1,-4,0,0)^T$, and $\mathcal{C}(2)=(1,1,-2,0,0,0, 0)^T$. The next columns of $V(G)$ are associated with $\mathcal{C}(i)$ for every $i$ that $T^G(i)=1$. Then, $V(G)$ is obtained as:
\begin{equation*}
V(G)=\begin{bmatrix}1 & 1 & 1 & 1 & 1 & 1 & 1\\ 1 & 1 & 1 & 1 & -1 & 1 & 1\\ 1 & 1 & 1 & -2 &  & -3 & 1\\1 & 1 & 1 &  &  &  & 1\\1 & 1 & -4 &  &  &  & 1\\1 & -5 &  &  &  &  & 1\\1 &  &  &  &  &  & -6\end{bmatrix}
\end{equation*}

\emph{Proof of Theorem \ref{th4}:} The proof follows by induction. First, note that for a threshold graph of size 1 which is an isolated node, $\Lambda(G)=(0)$ and $V(G)=[1]$. Now, consider $T^G$ for a threshold graph $G$ of size 2. Then, $T^G=(0,1)$ or $(0, 0)$. Then, it follows from Theorem \ref{th3} that $$V(G)=\begin{bmatrix}1 &1\\1 & -1\end{bmatrix},$$ which can be constructed by running Algorithm 1 as well. Now, assume that for any threshold graph $G'$ of size $n$, $V(G')$ can be obtained through Algorithm 1. Then, consider a threshold graph $G$ of size $n+1$. We want to prove that  $V(G)$ can be generated by running Algorithm 1. Let $T^G$ be the construction sequence of $G$. Then, we have either $T^G=(T^{G'},0)$ or $T^G=(T^{G'},1)$, where $T^{G'}$ is a construction sequence associated with a threshold graph $G'$ of size $n$. Thus, $V(G')$ can be provided by Algorithm 1. Let $\Lambda(G')=(0,\lambda_2,\ldots,\lambda_n)=(0,\Lambda')$, where $\Lambda'=(\lambda_2,\ldots,\lambda_n)$, and $0\leq  \lambda_2\leq \ldots \leq\lambda_n$. Moreover, let $V(G')=[\textbf{1}_n,V']$, where $V'=[\nu_2,\ldots,\nu_n]$. Now, first assume $T^G=(T^{G'},0)$. Thus, the node $n+1$ is added to the graph $G'$ through a union operation. Then, from Theorem \ref{th3}, $\Lambda(G)=(0,0,\Lambda')$, and $$V(G)=\begin{bmatrix}\textbf{1}_{n} & \textbf{1}_{n} & V'\\1 & -n & 0\end{bmatrix}.$$  Therefore, $V(G)$ can be constructed through Algorithm 1. In fact, according to this algorithm, since $T^G(n+1)=0$, $V_{:,2}(G)=\mathcal{C}(n)$, which is true. Now, assume that $T^G=(T^{G'},1)$ which means that the node $n+1$ is added to $G'$ through a join operation. According to Algorithm 1, since $T^G(n+1)=1$, $V_{:,n+1}(G)=\mathcal{C}(n)$. This can also be verified through Theorem \ref{th3} which implies that  $\Lambda(G)=(0,\Lambda',n+1)$, and $$V(G)=\begin{bmatrix}\textbf{1}_{n} & V' & \textbf{1}_{n} \\1 & 0 &-n \end{bmatrix}.$$       
\carre

We should note that for a threshold graph of size $n$, one can run Algorithm 1 in $\mathcal{O}(n)$.

\subsection{Controllability Analysis of Threshold Graphs}

We now consider a network with dynamics (\ref{e1}) with a connected threshold graph $G$. Furthermore, we assume that the input matrix $B$ is defined as  (\ref{B}). In particular, we proceed to characterize the minimal set of control nodes $V_C$ which renders the network controllable. 

Before presenting the control node selection method, let us introduce some more notation and present a lemma which is applied in the proof of the main result. For a connected threshold graph $G$ with the degree sequence $D(G)=(\tilde{d}_1^{p_1},\ldots,\tilde{d}_s^{p_s})$ and the spectrum $\Lambda(G)=(0,\tilde{\lambda}_1^{q_1},\ldots,\tilde{\lambda}_s^{q_s})$, let $V(G)=[\textbf{1}_n, \tilde{V}^{(1)},\ldots,\tilde{V}^{(s)}]$, where for every $1\leq i\leq s$, $\tilde{V}^{(i)}\in \mathbb{R}^{n\times q_i}$ is a matrix whose columns are the independent eigenvectors associated with the eigenvalue $\tilde{\lambda}_i$. Then, every vector $\tilde{\nu}_i=\tilde{V}^{(i)}C$, for some $C\in \mathbb{R}^{q_i} $, is an eigenvector  associated with $\tilde{\lambda}_i$. 

\begin{lem}
Consider a connected threshold graph with $D(G)=(\tilde{d}_1^{p_1},\ldots,\tilde{d}_s^{p_s})$ and $\Lambda(G)=(0,\tilde{\lambda}_1^{q_1},\ldots,\tilde{\lambda}_s^{q_s})$. If $s=2l$ (resp., $s=2l+1$), for some $l\in\mathbb{Z}$, then for $i\neq l$ (resp., $i\neq l+1$), there is some $1\leq k\leq n-1$ that \begin{equation}\tilde{V}^{(i)}=\begin{bmatrix}\mathcal{C}(k) & \ldots & \mathcal{C}(k+p_i-1) \end{bmatrix}. \label{e3}\end{equation} Moreover, if $i= l$ (resp., $i= l+1$), one can obtain that for some $1\leq k\leq n-1$,  \begin{equation}\tilde{V}^{(i)}=\begin{bmatrix}\mathcal{C}(k) & \ldots & \mathcal{C}(k+p_i-2)\end{bmatrix}.\label{e4} \end{equation}
\label{lem2}
\end{lem}

\emph{Proof:} We prove the result only for the case that $s=2l$ and $i\neq l$. The result for the other cases can be be proved in a similar way. Note that in this case, Theorem \ref{th2} implies  that $q_i=p_i$, that is, the multiplicity of the eigenvalue $\tilde{\lambda}_i$ is equal to the multiplicity of the $i$th degree. Then, $\tilde{V}^{(i)}\in \mathbb{R}^{n\times p_i}$. Moreover, from Lemma \ref{lem1}, $p_i$ nodes have the same degrees if they are successively indexed. In other words, if for the nodes $j_1, \ldots, j_{p_i}$, we have $d(j_1)=\ldots=d(j_{p_i})$, there is some $1\leq k\leq n-1$ that $j_r=k+r-1$, for $1\leq r\leq p_i$. Then, from the construction method of $V(G)$ in Algorithm 1, we have $\tilde{V}^{(i)}=\begin{bmatrix}\mathcal{C}(k) & \ldots & \mathcal{C}(k+p_i-1) \end{bmatrix}$. 
\carre

Now, let us partition the node set of a connected threshold graph into cells, such that the degrees of any two nodes in a cell are the same; while the degrees of two nodes from two different cells are different.  In the following, we show that the network is controllable if and only if from any cell, all nodes except one are chosen as control nodes. The procedure of the selection of the control nodes is presented as follows. 

\emph{Procedure 1}: Consider a connected threshold graph $G$ with the degree sequence $D(G)=(\tilde{d}_1^{p_1},\ldots, \tilde{d}_s^{p_s})$. For every $1\leq i\leq s$, let $K^{(i)}=\{1\leq j\leq n:d(j)=\tilde{d}_i\}$. Note that $|K^{(i)}|=p_i\geq 1$. Now, choose one node $k_i$ from every set $K^{(i)}$, $1\leq i\leq s$. Let $V'=\{k_1,\ldots,k_s\}$ and $V_C=V\setminus V'$. Then, $|V_C|=n-s$, where $n$ is the size of network, and $s$ is the number of distinct degrees of its nodes.

\begin{theo}
Consider a network with a connected threshold graph $G$ and dynamics (\ref{e1}) whose input matrix $B$ is described in (\ref{B}). Then, the network is controllable if $V_C$ is chosen through Procedure 1. Moreover, the minimum number of control nodes rendering the network controllable is $n-s$ which is also determined through the application of Procedure 1. 
\label{th5}
\end{theo}   

\emph{Proof:} Let $V_C=V\setminus V'$, and note that from Corollary \ref{corr1}, either $1\in V_C$ or $2\in V_C$. Now, assume that the network is not controllable. Then, for some $1\leq i\leq s$, there is a nonzero eigenvector $\tilde{\nu}_i$ associated with $\tilde{\lambda}_i$ such that $\tilde{\nu}_i^TB=\textbf{0}$. Then, for some nonzero $C\in\mathbb{R}^{{q}_i}$, one can write $\tilde{\nu}_i=\tilde{V}^{(i)}C$. Accordingly, from the PBH test, we should have $C^T(\tilde{V}^{(i)})^TB=\textbf{0}$. Note that $(\tilde{V}^{(i)})^TB=(\tilde{V}^{(i)}_{V_C,:})^T$, where $\tilde{V}^{(i)}_{V_C,:}$ is a submatrix of $\tilde{V}^{(i)}$ including its $j$th rows with all $j\in V_C$. From Lemma \ref{lem2}, for every $1\leq i\leq s$, $\tilde{V}^{(i)}$ in (\ref{e3}) has $p_i+1$ independent rows, that is, the rows $k,\ldots, k+p_i-1$ and one of the rows 1 and 2.  Then, if $V_C$ includes $p_i-1$ nodes from the $p_i$ nodes with the same degree along with one of nodes 1 and 2, then $\tilde{V}^{(i)}_{V_C,:}$ is full rank; thus, $C^T\tilde{V}^{(i)}_{V_C,:}=\textbf{0}$ implies that $C=0$; that is, $\tilde{\nu}_i=0$, which is a contradiction. For the second part of the theorem, we can do a similar argument and conclude that by choosing a set of control nodes with a size less than $n-s$, for some $1\leq i\leq s$,  $\tilde{V}^{(i)}_{V_C,:}$ is not full rank; then, $\tilde{\lambda}_i$ has some nonzero eigenvector $\tilde{\nu}_i$ that $\tilde{\nu}_i^TB=\textbf{0}$. Hence, the system would not be controllable.   
\carre
 
 As an example, consider the network with the threshold graph shown in Fig. \ref{ex1}. Applying Procedure 1, one can choose one of the nodes 1 and 2 and one of the nodes 5 and 6 as control nodes. For instance, we can have $B=[e_1,e_5]$. 
 
 \section{Controllability of Networks Defined on Cographs}
 
 In this section, we discuss the controllability of a network with dynamics (\ref{e1}) defined on a cograph. A cograph has a few definitions, all of which are  equivalent. Here, we describe a cograph by its associated \emph{cotree}. 
 
A cotree $\mathcal{T}$ associated with a cograph $G$ is a rooted tree whose leaves (i.e., the nodes with the degree one) correspond to the nodes of the cograph. Moreover, the internal nodes of a cotree (i.e., the nodes whose degree is bigger than one) are labeled with 0 or 1. Any subtree rooted at each node $z$ of $\mathcal{T}$ corresponds to an induced subgraph of $G$ defined on the leaves descending from $z$. If $z$ is a leaf of $\mathcal{T}$, the corresponding subgraph in $G$ is a graph of the single node $z$. In addition, to an internal node $z$ of $\mathcal{T}$ that is labeled 0, one can correspond a subgraph which is the union of subgraphs associated with the children of $z$. On the other hand, if $z$ is labeled 1, the corresponding subgraph is a join of subgraphs corresponding to  the children of $z$ \cite{corneil1985linear}.

We note that a cograph $G=(V,E)$ can be recognized in $\mathcal{O}(|V|+|E|)$, while its associated cotree can be constructed with similar computational efficiency \cite{corneil1985linear}.        

\begin{figure}[hbt]
\includegraphics[width=.43\textwidth]{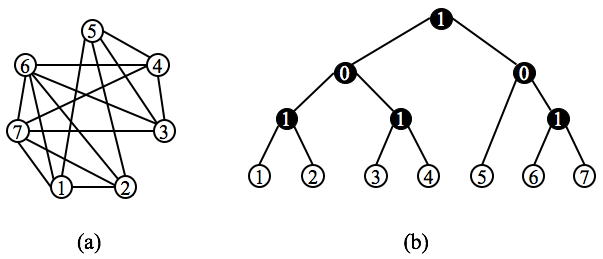}
\centering
\caption{a) A cograph $G$, b) Associated cotree $\mathcal{T}$.}
\label{cotree}
\end{figure}


In Fig. \ref{cotree}, an example of a cograph along with its associated cotree is illustrated. 

In order to characterize the eigenvalues and eigenvectors associated with a cograph, one can apply a  bottom-up tree computation on its associated cotree and by applying Theorem \ref{th3}, provide the spectrum of the cograph as well as its corresponding modal matrix in a polynomial time. By having the eigenspace of a network, its controllability problem can be addressed as follows.

Consider a connected cograph $G$. Let $\Lambda(G)=(\tilde{\lambda}_1^{q_1}, \ldots, \tilde{\lambda}_r^{q_r})$ be the spectrum of $G$ and $V(G)=[V^{(1)},\ldots, V^{(r)}]$ be its normalized modal matrix, where $\tilde{\lambda}_i$, $1\leq i\leq r$, is an eigenvalue of $L(G)$ with the multiplicity $q_i$. Moreover, $V^{(i)}\in \mathbb{R}^{n\times q_i}$, and $L(G)V^{(i)}=\tilde{\lambda}_iV^{(i)}$. Let $V^{(i)}=[V^{(i)}_1,\ldots,V^{(i)}_{q_i}]$, where $V^{(i)}_j\in \mathbb{R}^n$, for $1\leq j\leq q_i$. In addition, for some $1\leq k\leq r$, let $q_k$  be  the maximum multiplicity of eigenvalues of $G$. Now, for every $1\leq i\leq r$, add $q_k-q_i$ zero columns to each $V^{(i)}$ and define $\bar{V}^{i}=[V^{(i)}_{1},\ldots,V^{(i)}_{q_i}, \textbf{0}_{n\times (q_k-q_i)}]$. Then, we have the following condition for the controllability of the associated network.

\begin{theo}
A network with dynamics (\ref{e1}) which is defined on a cograph $G$ is controllable if $$B=[\sum_{i=1}^r \bar{V}_{:,1}^{(i)}, \ldots,\sum_{i=1}^r \bar{V}_{:,q_k}^{(i)}] \in \mathbb{R}^{n\times q_k }.$$ 
\label{th6}
\end{theo} 

 \emph{Proof:} Since $V(G)$ is normalized, $V(G)V^T(G)=V^T(G)V(G)=I_n$. Then, the equation $L(G)V(G)=V(G) \mbox{diag}(\Lambda(G))$ implies that $\mbox{diag}(\Lambda(G))=V^T(G)L(G)V(G)$. Let $D=\mbox{diag}(\Lambda(G))$. Since the controllability property is not influenced by the similarity transformation, the controllability of the pair $(L(G),B)$ is equivalent to the controllability of the pair $(D,\bar{B})$, where $\bar{B}=V^T(G)B$. Based on the PBH test, the pair $(D,\bar{B})$ is controllable if and only if for every $\tilde{\lambda}_i\in \Lambda(G)$, the matrix $[D-\tilde{\lambda}_iI,\bar{B}]$ is full rank. Accordingly, the pair $(D,\bar{B})$ is controllable if and only if the rows of $\bar{B}$ associated with the same diagonal entries of $D$ are independent. Now, for every $1\leq i\leq r$, let us define $E^{(i)}=[e_1,\ldots,e_{q_i}]$, where $e_j$, $1\leq j\leq q_i$,  is the $j$th column of $I_{q_k}$. We then choose $\bar{B}^T=[E^{(1)},\ldots,E^{(r)}]$. Thus, $B=V(G)\bar{B}=[\sum_{i=1}^r \bar{V}_{:,1}^{(i)}, \ldots,\sum_{i=1}^r \bar{V}_{:,q_k}^{(i)}]$.
 \carre

As an example, consider the cograph $G$ shown in Fig. \ref{cotree} (a). By having the associated cotree in Fig. \ref{cotree} (b) and applying Theorem \ref{th3}, one can obtain $\Lambda(G)=(0,3,4,5^2,6,7)$. Moreover, we have:

\begin{align*}
V(G)=\begin{bmatrix} 
1 & 2 & & & 1 & & 3\\
1 & 2 & & & -1 & & 3 \\
1 & -2 & & 1 &  & & 3\\
1 & -2 & & -1 &  & & 3\\
1 &  &2 & &  & & -4\\
1 & & -1& &  & 1& -4\\
1 &  & -1 & & & -1& -4
\end{bmatrix}.
\end{align*}

Accordingly,  $V^{(1)}=V_{:,1}(G)$, $V^{(2)}=V_{:,2}(G)$, $V^{(3)}=V_{:,3}(G)$, $V^{(4)}=[V_{:,4}(G),V_{:,5}(G)]$, $V^{(5)}=V_{:,6}(G)$, and $V^{(6)}=V_{:,7}(G)$. Then from Theorem \ref{th6}, one can obtain $B=[\sum_{i=1}^6 V_{:,1}^{(i)}, V_{:,2}^{(4)}]$. Thus, we have:

\begin{align*}
B^T=\begin{bmatrix} 
6 & 6 & 3 & 1 & -1 & -3 & -5 \\ 1  &-1 & 0& 0&0 &0 &0 
\end{bmatrix}.
\end{align*}

Note that the entries of the input matrix $B$ obtained from Theorem \ref{th6} are all integer. 


\section{Conclusion}

In the first part of this paper, the controllability of an LTI  network with the Laplacian dynamics on a threshold graph has been examined.  In this direction, an efficient algorithm for characterizing the modal matrix associated with the Laplacian of a threshold graph has been presented. Subsequently, by assuming that any input signal can be injected into one node only, necessary and sufficient conditions for the controllability of this class of networks has been established. Furthermore, the paper examined  the controllability problem of general cographs; it is shown that an input matrix with the minimum rank that renders the network controllable can be found in a  polynomial time.    

\bibliographystyle{IEEEtran}
\bibliography{library}

\end{document}